\documentclass[11pt]{article}
\usepackage{amsmath,amssymb,amsthm}

\def\g{\mathfrak{g}}
\def\C{\mathbb{C}}

\def\Z{\mathbb{Z}}
\def\N{\mathbb{N}}
\def\qed{$\hfill \square$}
\def\l{\mathfrak{l}}

\def\m{\mathfrak{m}}
\def\zom{\frac{1}{4}(Z_{\Omega}-Z_{\Omega_\l})}
\def\h{\mathfrak{h}}

\newtheorem{theo}{Theorem}
\newtheorem{theoA}{Theorem A}
\newtheorem{prop}{Proposition}[section]
\newtheorem{propA}{Proposition A}
\newtheorem{lemA}{Lemma A}
\newtheorem{lem}{Lemma}[section]
\newtheorem{cor}{Corollary}[section]

\numberwithin{equation}{section}

\begin{document}

\title{On the Moduli Space of Classical Dynamical r-matrices}

\author{Pavel Etingof and Olivier Schiffmann}
\date{}
\maketitle

\paragraph{Introduction.} A classical dynamical r-matrix is an $\l$-equivariant
function 
$r:\;\l^*\to \g\otimes \g$ (where $\l$, $\g$ are Lie algebras),
such that $r^{21}+r=\Omega$ is $\g$-invariant, which satisfies 
the classical dynamical Yang-Baxter equation (CDYBE). 
CDYBE is a differential equation, 
which generalizes the usual classical Yang-Baxter equation. 
It was introduced in 1994 by G.Felder \cite{Fe}, in the context of 
conformal field theory. Solutions of CDYBE 
and their quantizations appear naturally in 
several mathematical theories: the theory of 
integrable systems, special functions, representation theory
(see \cite{ES} for a review). 

Since classical dynamical r-matrices were introduced, several authors tried to 
study and classify them (\cite{EV},\cite{S},\cite{Xu}).  
The goal of this paper is to describe the moduli space of 
classical dynamical r-matrices modulo gauge transformations. 
In particular, we improve and generalize the results 
of \cite{EV}, \cite{S}, as well as correct some errors that 
occurred in these papers (See remarks 3 and 5).

The main achievement of this paper, compared to the previous ones, 
is that its results are valid for dynamical r-matrices
for a {\bf nonabelian} Lie algebra $\l$. 
It turns out that this generalization not only brings in 
new interesting examples (see \cite{EV},\cite{AM}) but also
makes the general theory much more clear and natural. 

The composition of the paper is as follows. 

In Section 1, we recall the definition of a dynamical r-matrix. 
 
In Section 2, we extend to the nonabelian case the notion of 
a gauge transformation of dynamical r-matrices, introduced in \cite{EV}. 

In Section 3, we decribe the space of dynamical r-matrices modulo 
gauge transformations (the moduli space).
Here we formulate our main theorem, stating that 
under some technical conditions, 
the moduli space can be identified with a certain explicitly 
given affine variety. 
For instance, if $\l=\g$, this variety 
consists of one point, which is the Alekseev-Meinrenken solution 
\cite{AM} (for semisimple Lie algebras, it was also constructed in \cite{EV}).

In Section 4 we prove the main theorem. 

In the appendix, we construct a generalization of the Alekseev-Meinrenken
classical dynamical r-matrix, associated to any finite-dimensional Lie
algebra $\g$ with a nondegenerate invariant form and an automorphism $B$
of $\g$ of finite order which preserves this form. 

\section{The Dynamical Yang-Baxter Equation}
\paragraph{}Let $\g$ be a Lie algebra over 
$\C$, and $\l \subset \g$ a finite dimensional Lie 
subalgebra. Let $x_1,...,x_r$ be a basis of $\l$. 

Let $D \subset
\l^*$ be the formal neighborhood of $0$.
Let $V$ be a complex vector space. By functions from $D$ to $V$ 
we will mean elements of the space $V[[x_1,...,x_r]]$, where we regard
$x_i$ as coordinates on $D$. Finally, if $\omega \in \Omega^k(D,V)$ is a
$k$-form with values in any vector space $V$, we denote by $\overline{\omega}:
D \to \Lambda^k \l \otimes V$ the associated function. For an element
 $r \in \g \otimes \g$ we define the classical Yang-Baxter operator
$$CYB(r)= [r^{12},r^{13}]+[r^{12},r^{23}]+[r^{13},r^{23}].$$
The 
\textbf{classical dynamical Yang-Baxter equation}
(CDYBE) is the following differential
equation for 
an $\l$-equivariant function
$r: D \to \g \otimes
\g$~:
\begin{equation}\label{E:06}
\mathrm{Alt}(\overline{dr}) + CYB(r)=0,
\end{equation}
where 
for $x \in \g^{\otimes 3}$, we let $\mathrm{Alt}(x)=x^{123} - x^{213}
+x^{312}$.

It is useful to consider solutions of 
CDYBE which satisfy an additional quasi-unitarity condition:
\begin{equation}
 r + r^{21}=\Omega \in
(S^2\g)^{\g}. \label{E:08}
\end{equation}
It is easy to show that if $r$ satisfies CDYBE and 
the quasi-unitarity condition then $\Omega$ is a constant
function of $\lambda$. 

An $\l$-equivariant solution of CDYBE 
which satisfies the quasi-unitarity condition
is called a {\it dynamical r-matrix}.
The set of all dynamical r-matrices satisfying (\ref{E:08})
will be denoted by $Dynr(\g,\l,\Omega)$.

\section{Gauge transformations}
\paragraph{}Here we will reproduce some results from \cite{EV}, 
but unlike \cite{EV}, we will not assume that $\l$ is
 abelian. We will assume, however, that $\g$ is finite
 dimensional.  

Let $G$ be the simply connected complex Lie group such that  
$\mathrm{Lie}(G)=\g$. Let $g: D \to G$ 
be any regular, $\l$-equivariant map. Consider the $1$-form
$\eta_g=g^{-1}dg$ and set
$\zeta_g= [\overline{\eta_g}^{12},\overline{\eta_g}^{13}]$. Define an
$\l$-equivariant function
$\tau_g: D \to \Lambda^2\g$ by the formula
$\tau_g(\lambda)
=(\lambda \otimes 1 \otimes 1)\zeta_g(\lambda)$.
For any 
$\l$-equivariant function $r: D \to
\g\otimes
\g$ we set
\begin{equation}\label{E:un}
r^g=(g \otimes g)(r-\overline{\eta_g}+\overline{\eta_g}^{21}+\tau_g)
(g^{-1} \otimes g^{-1}).
\end{equation}
The following theorem is a nonabelian generalization of 
Proposition 1.2 of \cite{EV}. 

\begin{prop} 
The function $r$ is a dynamical r-matrix if and only if the
function $r^g$ is.
\end{prop}
\noindent
\textit{Proof.} Let us show that if $r$ is a dynamical r-matrix then so is
$r^g$. The other direction is analogous.
Let $X =(D\times
G \times D,\{\,,\,\})$ be the dynamical Poisson groupoid associated to $r$ in
\cite{EV}. Consider the automorphism $\sigma$ of $X$ given by
$\sigma(u_1,x,u_2)=(u_1,g(u_1)xg(u_2)^{-1},u_2)$. Then $\sigma$ transforms
$\{\,,\,\}$ into the Poisson bracket $\{f,g\}_{\sigma}=\sigma^{-1}\{\sigma f,
\sigma g\}$. It is straightforward to calculate that the corresponding transformation at the
level of dynamical r-matrices is exactly (\ref{E:un}).\qed

\paragraph{}The transformation 
$r\to r^g$ is called a gauge transformation.
Note that (\ref{E:un}) defines an action of the group $\mathrm{Map}(D,G)^\l$
on $Dynr(\g,\l,\Omega)$, i.e we have $(r^{g_1})^{g_2}=r^{g_2g_1}$ for any
$g_1,g_2 \in \mathrm{Map}(D,G)^\l$ and $r \in Dynr(\g,\l,\Omega)$.
Let us denote by
$\mathrm{Map}_0(D,G)^\l$ the subgroup consisting of maps $g$ satisfying
$g(0)=1$. We would like to understand the moduli
space 
$$\mathcal{M}(\g,\l,\Omega)=Dynr(\g,\l,\Omega)/\mathrm{Map}_0(D,G)^\l.$$
In the triangular case (i.e when $\Omega=0$) this space was considered by
P. Xu in \cite{Xu}.
\paragraph{Remark 1.} It is clear that $\mathrm{Map}(D,G)^\l/
\mathrm{Map}_0(D,G)^\l\simeq G^\l$. Hence the
complete moduli space $\overline{\mathcal{M}}(\g,\l,
\Omega)=Dynr(\g,\l,\Omega)/\mathrm{Map}(D,G)^\l$ is equal to
$\mathcal{M}(\g,\l,\Omega)/G^\l$ where $g \in G^\l$ acts by
$r^g=\mathrm{Ad}(g \otimes g)(r)$. 

\section{The structure of $\mathcal{M}(\g,\l,\Omega)$}

From now on we will assume that 
\begin{enumerate}
\item[i)]
$ \l\subset \g\; \mathrm{has \;an\;} 
\l\mathrm{-invariant\;complement}\;\m.$
\end{enumerate}
The following 
theorem is a generalization 
of Theorem 1.4 in \cite{EV}. 
It shows that the space of dynamical r-matrices is, up
to gauge equivalence, finite dimensional. 

\begin{theo} Let $\rho,r:D \to \g^{\otimes 2}$ be two dynamical
r-matrices such that $r(0)=\rho(0)$. Then there
exists $g \in \mathrm{Map}(D,G)^\l$ such that $\rho=r^g$.
\end{theo}
\noindent
The proof is a generalization of the proof in \cite{EV}.
Before giving it we state the
following auxiliary result.
\begin{lem}[equivariant Poincar\'e lemma] Let $\l$ be a finite-dimensional Lie
algebra, $V$ a finite-dimensional
$\l$-module, $k \geq 1$ and
$\omega \in \Omega^k(D,V)$ an $\l$-equivariant closed $k$-form with values
in $V$. Then there exists an $\l$-equivariant $k-1$-form $\zeta \in
\Omega^{k-1}(D,V)$ such that $d\zeta=\omega$.\end{lem}
\noindent
\textit{Proof.} The proof is the same as that for the usual Poincar\'e lemma.
It is enough to assume that $\omega$ is homogeneous, of degree
$l \in \N$. Let $E=\sum_i x_i \frac{\partial}{\partial x_i}$ be the Euler
vector field on $D$. Then by Cartan's homotopy formula,
$$l\omega=L_E \omega=i_E d\omega+d i_E \omega=d(i_E \omega)$$
and we can set $\zeta=i_E\omega/l$. Note that $E$ is $\l$-equivariant, hence
so is $\zeta$. \qed\\

\vspace{1 mm}

\noindent
\textit{Proof of Theorem 1.} The dynamical r-matrices $r,\rho$ are by
definition formal
power series in the variables $x_i$. Let us assume that the statement
of the theorem
holds modulo terms of degree $\geq K$. Let $g_k: U \to G$ be a
gauge transformation such that $E:=r^{g_k}-\rho$ has degree $\geq K$ and
let $E_K$ be the homogeneous component of $E$ of degree $K$. Then
$\mathrm{Alt}\;(\overline{dE_K})=\big[ CYB(r^{g_K})-CYB(\rho)\big]_{K-1}$
where $\big[\cdot\big]_{K-1}$ denotes the homogeneous component of degree
$K-1$. But
$\big[r^{g_K}-\rho\big]_{l}=0$ for all $l <K$ by assumption, hence 
\begin{equation}\label{E:t1}
\mathrm{Alt}\;(\overline{dE_K})=0.
\end{equation}
\begin{lem} The exists an $\l$-equivariant closed 1-form $\zeta \in
\Omega^1(D,\g)$ such that $E_K=\overline{\zeta}^{21}-\overline{\zeta}$.
\end{lem}
\noindent
\textit{Proof.} Let us write $E_K=E_{\l\l}+E_{\l\m}-E_{\l\m}^{21}+E_{\m\m}$
where $E_{\l\l} \in \Lambda^2\l$, $E_{\l\m} \in \l \otimes \m$ and
$E_{\m\m}\in \Lambda^2\m$. From (\ref{E:t1}) it follows that $dE_{\m\m}=0$
hence $E_{\m\m}=0$. Now let $\xi \in \Omega^1(D,\m)$ be such that
$\overline{\xi}=E_{\l\m}$. Then (\ref{E:t1}) implies that $\xi$ is closed.
Note that the assumption i) guarantees that $\xi$ is equivariant.
Finally, let $\omega \in \Omega^2(D,\C)$ be such that $\overline{\omega}
=E_{\l\l}$. Then (\ref{E:t1}) says that $\omega$ is closed. By the equivariant
Poincar\'e lemma, there exists an equivariant 1-form $\eta$ such that
$d\eta=\omega$. Set $\theta=d\overline{\eta}$, so that $\overline{\theta}
-\overline{\theta}^{21}=\overline{\omega}$. Then $\zeta=\xi + \theta$ satisfies
the conditions of the lemma.\qed
\paragraph{}We now conclude the proof of Theorem 1. Let $\chi: D \to \g$ be any
$\l$-equivariant function of order $K+1$ such that $d\chi=\zeta$. Set
$g=e^\chi$. Then $\eta_g=g^{-1}dg$ is of order $\geq K$ and 
$\overline{\zeta}-\overline{\eta_g}$ is of order $\geq K+1$.
But then $\overline{\zeta}^{21}-\overline{\zeta}-(\overline{\eta_g}^{21}-
\overline{\eta_g}+\tau_g)$ is also of order $\geq K+1$.
Set $g_{K+1}=gg_K$. Then, by the above $r^{g_{K+1}}-\rho$ is of degree
$\geq K+1$. The proof follows by induction.\qed

\begin{prop}Any dynamical r-matrix $r$ is gauge-equivalent to a dynamical
r-matrix $\rho$ such that $\rho(0) \in \frac{\Omega}{2} + \Lambda^2\m$.
\end{prop}
\noindent
\textit{Proof.} Let $\overline{\eta}_0 \in \l \otimes \g$ such that
$r(0)-\overline{\eta}_0+\overline{\eta}_0^{21} \in \frac{\Omega}{2} +
\Lambda^2\m$. Since $\m$ is $\l$-invariant, we have $\overline{\eta}_0
\in (\l \otimes \g)^\l$. By the equivariant Poincar\'e lemma, there exists
an equivariant function $\chi: D \to \g$ satisfying $\chi(0)=0$, $d\chi=
\eta_0$. Set $g=e^\chi$. Then $\rho:=r^g$ satisfies the CDYBE and $\rho(0)
\in \frac{\Omega}{2} + \Lambda^2\m$. \qed

\paragraph{}Consider the following algebraic variety
$$ \mathcal{M}_{\Omega}=
\{x \in \frac{\Omega}{2}+(\Lambda^2\m)^\l\;|
CYB(x)=0\;in\;\Lambda^3(\g/\l)\}.$$
It is immediate from (\ref{E:un}) that if $\rho$ and $r$ are 
gauge-equivalent and if $\rho(0) \in \frac{\Omega}{2} + \Lambda^2\m$ and
$ r(0)\in\frac{\Omega}{2} + \Lambda^2\m$ then $r(0)=\rho(0)$. Moreover,
it follows from the CDYBE (\ref{E:06}) that for every dynamical r-matrix
$r \in Dynr(\g,\l,\Omega)$ such that $r(0) \in \frac{\Omega}{2} + \Lambda^2\m$
we have $r(0) \in \mathcal{M}_{\Omega}$.

Hence Theorem 1 and Proposition 3.1 give the following corollary.
\begin{cor}The map $\mathcal{M}(\g,\l,\Omega) \to \mathcal{M}_{\Omega}$
which sends a class
$\mathcal{C}$ to $r(0)$ where $r \in 
\mathcal{C}$ is any representative such that $r(0) \in
\frac{\Omega}{2} + \Lambda^2\m$, is an embedding.
\end{cor}

\paragraph{Remark 2.} If condition i) fails then the space
$\mathcal{M}(\g,\l,\Omega)$ may be
infinite-dimensional. This is demonstrated by the following
example due to P. Xu \cite{Xu}. Let $\g=\C x \oplus \C y$ be the
two-dimensional
Lie algebra with $[x,y]=y$, and set $\l=\C y$. Then $\Lambda^3 \g=0$ and
$\Lambda^2\g$ is a trivial $\l$-module. Thus
any function $r: D \to \Lambda^2 \g$ is a dynamical r-matrix. On the other
hand, $\g^\l=\l$ and all gauge transformations act trivially.

\paragraph{Remark 3.} We would like to use this opportunity to
correct the statement of Theorem 1.4 of \cite{EV}. 
This theorem is incorrect as stated (as shown by Xu's counterexample, see
Remark 1).
The mistake is in the proof of Lemma 1.5, which
uses the incorrect statement that 
\begin{equation}\label{E:mistake}
(\g\otimes \l\oplus \l\otimes \g)^\l=
(\g^\l\otimes \l\oplus \l\otimes \g^{\l})^\l
\end{equation}
for commutative $\l$. This statement, however, is correct 
with the additional assumption i);
in this case Theorem 1.4 of \cite{EV} and its proof are correct,
and Theorem 1.4 of \cite{EV} is a special case of 
Theorem 1 above. 

\paragraph{}Now suppose that $\l=\g$. Note that i) automatically holds in this
case. Then by Proposition 3.1 and Theorem 1 there is at most one
gauge-equivalence class of dynamical r-matrices $r: D \to \frac{\Omega}{2}
+ \Lambda^2 \g$. Such a class in fact always exists, as was discovered by
Alekseev and Meinrenken \cite{AM}. A representative of this class is
constructed as follows.
\paragraph{}Let $\g_\Omega$ be the ideal of $\g$
spanned by the components of $\Omega$, and let $D_\Omega$ be the formal
neighborhood of $0$ in $\g^*_\Omega$. Let us identify $\g_\Omega$ with
$\g^*_\Omega$ via $\Omega$. Set $f(s)=\frac{1}{s}-\frac{1}{2}
\mathrm{cotanh}(\frac{s}{2})$. Then $f$ is smooth
at the origin.
Consider the following map
\begin{align*}
T:\;D_\Omega &\to \mathrm{End}(\g) \simeq \g^* \otimes \g \simeq
 \g \otimes \g\\
u &\mapsto f(\mathrm{ad}\;\mu)
\end{align*}
Let $\pi^*: \g^* \to \g_\Omega^*$
be the projection and set
$$r^\g_{AM}=\frac{\Omega}{2} + T \circ \pi^* : D \to \g \otimes \g.$$
\begin{theo}[\cite{AM}] The map $r^\g_{AM}$ is a dynamical r-matrix.
\end{theo}
This theorem is proved in \cite{AM} in the case of compact Lie algebras, but
the proof can be adapted to the general case. Another proof is given in the
appendix.
\begin{cor} The moduli space $\mathcal{M}(\g,\g,\Omega)$ consists of the single
class $r^\g_{AM}$.\end{cor}
\paragraph{Remark 4.} When $\g$ is a simple Lie algebra and $\l=\g$ these
results easily follow from \cite{EV}, Section 3.8. 
\paragraph{}We will now show that, under some technical conditions
on $\Omega$, the embedding defined in Corollary 3.1 is actually an isomorphism.
From now on we assume that
\begin{enumerate}
\item[ii)] We have $\Omega \in (\l \otimes \l) \oplus (\m \otimes \m)$.
\end{enumerate}
We will write $\Omega_\l$ (resp. $\Omega_{\m}$) for the component of
$\Omega$.\\
\hbox to1em{\hfill}
Condition ii) is satisfied in particular in the triangular case
($\Omega=0$). It is
also satisfied when $\l$ acts semisimply on $\g$ and when the restriction of
the inverse form $(\,,\,)=\Omega^{-1}$ to $\l_\Omega=\l\cap\g_\Omega$ is
nondegenerate. Indeed, let $\g'$ be an $\l$-invariant complement of
$\l + \g_\Omega$ in $\g$ and let $\m_\Omega$ be the orthogonal complement
of $\l_\Omega$ i $\g_\Omega$. Then $\m=\g'\oplus\m_\Omega$ satisfies
conditions i) and ii).

\begin{prop}Any dynamical r-matrix $r$ is gauge-equivalent to a dynamical
r-matrix of the form $\rho=r^\l_{AM} + \frac{\Omega_{\m}}{2} + t$ with
$t : D \to \Lambda^2\m$.\end{prop}
\noindent
\textit{Proof.} By Proposition 3.1 there exists a dynamical r-matrix $\rho_0$
gauge-equivalent to $r$ such that $\rho_0(0) \in \frac{\Omega}{2}
+\Lambda^2\m$. We will first construct a sequence of gauge transformations
$g_i,\;i=1,\ldots$ such that $\rho_0^{g_i} \in \frac{\Omega}{2}+
\big(\Lambda^2\l \oplus \Lambda^2\m\big)$ modulo terms of degree $\geq i$.
We set $g_1=1$. Suppose that we have constructed $g_i$ and let $E_i$ be the
term of degree exactly $i$ of $\rho_0^{g_i}$. From the CDYBE we have
\begin{equation}\label{E:deux}
-\mathrm{Alt}\;(\overline{dE_i})=\big[CYB(\rho_0^{g_i})\big]_{i-1}
\end{equation}
where $[\cdot]_{i-1}$ denotes the component of degree $i-1$. But by our
assumption we have $\rho_0^{g_i} \in \frac{\Omega}{2} + \big( \Lambda^2\l 
\oplus \Lambda^2\m
\big)$ in degrees $\leq i-1$. Using the $\l$-invariance of $\m$ it is easy to
see that this implies that
\begin{equation}\label{E:trois}
\big[CYB(\rho_0^{g_i})\big]_{i-1} \in \mathrm{Alt}\;\big( (\l \otimes \m 
\otimes \m)
\oplus (\l \otimes \l \otimes \l) \oplus (\m \otimes \m\otimes\m)\big).
\end{equation}
Let $\xi \in \Omega^1(D,\m)$ such that $E_i + \overline{\xi}^{21}
-\overline{\xi} \in \Lambda^2\l \oplus \Lambda^2\m$. Then from (\ref{E:deux})
and (\ref{E:trois}) it follows that $d\xi=0$. By the equivariant Poincar\'e
lemma there exists an equivariant map $\chi : D \to \m$ such that $\xi=d\chi$.
Moreover, $\xi$ is of degree $\geq i$, hence $\chi$ is of degree $\geq i+1$.
Now set $g=e^\chi$. Then $\eta_g-\xi$ is of order $\geq i+1$. Thus
$$(g \otimes g)\big(\rho_0^{g_i}+\overline{\eta_g}^{21}-\overline{\eta_g}
+\tau_g\big)(g^{-1} \otimes g^{-1})$$
is in $\frac{\Omega}{2}+\big(\Lambda^2\l \oplus \Lambda^2\m\big)$
modulo terms of degree $\geq i+1$, and we put $g_{i+1}=g_ig$. This allows to
define the sequence $g_i$ inductively.\\
\hbox to1em{\hfill}It is clear that the sequence $\rho_0^{g_i}$ converges,
in the sense of formal power series, to a dynamical r-matrix $\rho_1$ which is
gauge-equivalent to $\rho_0$. Moreover $\rho_1$ takes values in
$\frac{\Omega}{2}+
\big(\Lambda^2\l \oplus \Lambda^2\m\big)$ by construction. Let us write
$\rho_1=\rho_1^\l + \rho_1^{\m}$ where $\rho_1^\l$ and $\rho_1^{\m}$
take values in $\l \otimes \l$ and $\m \otimes \m$ respectively.
Observe that $\rho_1^\l: D \to 
\frac{\Omega_\l}{2} + \Lambda^2\l$ is itself a dynamical r-matrix. Hence
by Corollary 3.2 we can
perform a gauge-transformation for $\l$ to
reduce it to $r^\l_{AM}$.\qed
\paragraph{}The following theorem is a generalization to the nonabelian case
of \cite{S}, Theorem 3, and will be proved in the next section.
\begin{theo} Let $r_0 \in \mathcal{M}_\Omega$. Then there exists a unique
dynamical r-matrix $r=r^\l_{AM} + \frac{\Omega_\m}{2}+t$ with
$t: D \to \Lambda^2\m$, such that $r(0)=r_0$.
\end{theo}
\begin{cor} Under conditions i) and ii) the moduli
space $\mathcal{M}(\g,\l,\Omega)$ of gauge-equivalence classes of dynamical
r-matrices is isomorphic to $\mathcal{M}_{\Omega}$.
\end{cor}

\paragraph{Remark 5.} We use this opportunity to correct
the statement of Theorem 3 in \cite{S} which is false as stated. The mistake
is in the proof of Lemma 1, which uses the incorrect statement
(\ref{E:mistake}). However, the theorem and its proof are correct
if one makes 
in addition the assumption i). In this case it is a special case of Theorem 3
above.
Moreover the genericity assumption made in \cite{S} Theorem 3 is not necessary,
as the flow constructed in \cite{S} Lemma 2 is well-defined on the whole 
$(\Lambda^2\m)^\l$.
\paragraph{Remark 6.} Let us identify $\m$ with $\g/\l$ via the decomposition
$\g=\l \oplus \m$. This allows to define an action of $G^\l$ on $\m$, hence
also an action of $G^\l$ on $\mathcal{M}_\Omega$.
It is clear from (\ref{E:un}) that the
isomorphism $\mathcal{M}(\g,\l,\Omega) \simeq \mathcal{M}_\Omega$ is
$G^\l$-equivariant. In particular, $\overline{\mathcal{M}}(\g,\l,\Omega)
\simeq \mathcal{M}_\Omega/G^\l$.

\section{Proof of Theorem 3}
\textit{Proof of Theorem 3.}
We will construct by induction a formal power series
$t=\sum_k t_k$ with $t_k : D \to \Lambda^2\m$ of degree $k$,
such that $r=r^{\l}_{AM}+\frac{\Omega_\m}{2}+t$ is a dynamical
r-matrix satisfying $r(0)=r_0$. Set $t_0=r_0-\frac{\Omega}{2} \in \Lambda^2\m$
and let us suppose that we have defined an $\l$-equivariant
polynomial $t_{<k}=\sum_{l<k} t_l$. Set $s=r^\l_{AM}-\frac{\Omega_\l}{2}$,
$Z_\Omega=CYB(\Omega)$ and $Z_{\Omega_\l}=CYB(\Omega_\l)$.
Then the CDYBE for $r^\l_{AM}$ is equivalent to the following equation for
$s$ :
\begin{equation}\label{E:mcdybe}
\mathrm{Alt}\;(\overline{ds})+CYB(s)+\frac{1}{4}Z_{\Omega_\l}=0.
\end{equation}
Let $\pi : \g \to \l$ be the projection along $\m$.
Consider, for $l \leq k$ the following system of differential equations for
$i=1,\ldots r$.
\begin{equation}\label{E:sept}
\frac{\partial t_l}{\partial x_i^*}=-(x_i^* \otimes 1 \otimes 1)
\bigg[[t_{<l}^{12},t_{<l}^{13}]+[s^{12}+s^{13},t_{<l}]+
\frac{1}{4}(Z_\Omega-Z_{\Omega_\l})\bigg]_{l-1}
\tag{$E_l$}
\end{equation}
where by definition $x^*(y)=x^*(\pi(y))$ for all $x^*\in \l^*$, $y \in \g$.
\begin{lem} Suppose that (\ref{E:sept}) is satisfied for all $l< k$. Then
($E_k$) admits a unique solution $t_k$ of degree $k$, which
is $\l$-equivariant.\end{lem}
\noindent
\textit{Proof.} By the equivariant Poincar\'e lemma, it is enough to show that
\begin{equation}\label{E:dix}
\begin{split}
\frac{\partial}{\partial x^*_j}&(x_i^* \otimes 1 \otimes 1)
\bigg\{[t_{<k}^{12},t_{<k}^{13}]+[s^{12}+s^{13},t_{<k}]+
\frac{1}{4}(Z_{\Omega}-Z_{\Omega_\l})\bigg\}\\
&=\frac{\partial}{\partial x_i^*}(x_j^* \otimes 1 \otimes 1)
\bigg\{[t_{<k}^{12},t_{<k}^{13}]+[s^{12}+s^{13},t_{<k}]+
\frac{1}{4}(Z_\Omega-Z_{\Omega_\l})\bigg\}
\end{split}
\end{equation}
Let us write $\partial_i$ for $\frac{\partial}{\partial x^*_i}$ and
$t$ for $t_{<k}$.
All equations below will be understood modulo
terms of degree $\geq k$. Let $X_i$ and $X_j$ denote the r.h.s and l.h.s
of (\ref{E:dix}).
Using the assumption that $t$ is a solution of the system (\ref{E:sept}) for
all $l <k$, we have
\begin{equation}\label{E:G}
\begin{split}
X_i&-X_j=\\
&=(x_i^*\otimes x_j^* \otimes 1\otimes 1)\bigg\{
\big[ [t^{12},t^{13}] +[s^{12}+s^{13},t^{23}]+\zom^{123},t^{24}\big]\\
&\qquad +\big[t^{23},[t^{12},t^{14}]+[s^{12}+s^{14},t^{24}] +\zom^{124}\big]\\
&\qquad +\big[ s^{23}+s^{24},[t^{13},t^{14}]+[s^{13}+s^{14},t^{34}]+
\zom^{134}\big]\\
&\qquad-[\partial_i(s^{23}+s^{24}),t^{34}]+[\partial_j(s^{13}+s^{14}),t^{34}]\\
&\qquad-\big[ -[t^{12},t^{23}]
+[-s^{12}+s^{23},t^{13}]-\zom^{123},t^{14}\big]\\
&\qquad-\big[t^{13},-[t^{12},t^{24}]+[-s^{12}+s^{24},t^{14}] -\zom^{124}\big]\\
&\qquad-\big[ s^{13}+s^{14},[t^{23},t^{24}]+[s^{23}+s^{24},t^{34}]
+\zom^{234}\big]
\bigg\}.
\end{split}
\end{equation}
By the Jacobi identity we have
\begin{equation}\label{E:A}
\big[ [t^{12},t^{13}],t^{24}\big] +\big[ t^{13},[t^{12},t^{24}]\big]=
\big[ [t^{12},t^{23}],t^{14}\big] +\big[ t^{23},[t^{12},t^{14}]\big]=0.
\end{equation}
Moreover,
\begin{equation}\label{E:B}
\begin{split}
(x_i^*\otimes x_j^* \otimes 1\otimes 1)\bigg\{
&[\frac{1}{4}(Z_{\Omega_\l})^{123},t^{24}]+[t^{23},\frac{1}{4}
(Z_{\Omega_{\l}})^{124}]\\
&+[\frac{1}{4}(Z_{\Omega_{\l}})^{123},t^{14}]+[t^{13},\frac{1}{4}
(Z_{\Omega_{\l}})^{124}]
\bigg\}
=0
\end{split}
\end{equation}
since $Z_{\Omega_\l} \in \Lambda^3\l$, $t \in \Lambda^2\m$ and $\m$ is
$\l$-invariant. Furthermore, $Z_{\Omega}$ is \\
$\g$-invariant, hence
\begin{align}\label{E:C}
(x_i^* \otimes x_j^* \otimes 1 \otimes 1)\bigg\{
&[\frac{1}{4}(Z_{\Omega})^{123},t^{24}]+[\frac{1}{4}
(Z_{\Omega})^{123},t^{14}]\notag\\
&+[t^{23},\frac{1}{4}(Z_{\Omega})^{124}]+[t^{13},\frac{1}{4}
(Z_{\Omega})^{124}]\bigg\}\notag\\
=(x_i^* \otimes x_j^* \otimes 1 \otimes& 1)\bigg\{
-[\frac{1}{4}(Z_{\Omega})^{123},t^{34}]
+[t^{34},\frac{1}{4}(Z_{\Omega})^{124}]\bigg\}\notag\\
=-(x_i^* \otimes x_j^* \otimes 1 \otimes& 1)\bigg\{
[\frac{1}{4}\big((Z_{\Omega_\l})^{123}+(Z_{\Omega_\l})^{124}\big),t^{34}]
\bigg\}.
\end{align}
In a similar way, $Z_{\Omega}$ and $Z_{\Omega_{\l}}$ are $\l$-invariant, hence
\begin{align}\label{E:D}
[s^{23}+s^{24},\frac{1}{4}(Z_{\Omega}-Z_{\Omega_{\l}})^{134}]&=
[s^{12},\frac{1}{4}(Z_{\Omega}-Z_{\Omega_{\l}})^{134}],\notag\\
-[s^{13}+s^{14},\frac{1}{4}(Z_{\Omega}-Z_{\Omega_{\l}})^{234}]&=
[s^{12},\frac{1}{4}(Z_{\Omega}-Z_{\Omega_{\l}})^{234}].
\end{align}
From the Jacobi identity again we deduce
\begin{align}\label{E:E}
\big[[s^{12},t^{23}],t^{24}\big] +\big[t^{23},[s^{12},t^{24}]\big]&
=\big[s^{12},[t^{23},t^{24}]\big],\notag\\
\big[[s^{12},t^{13}],t^{14}\big] +\big[t^{13},[s^{12},t^{14}]\big]&
=\big[s^{12},[t^{13},t^{14}]\big],\notag\\
\big[[s^{13},t^{23}],t^{24}\big]-\big[s^{13},[t^{23},t^{24}]\big]&=0,\notag\\
\big[s^{24},[t^{13},t^{14}]\big]-\big[t^{13},[s^{24},t^{14}]\big]&=0,\notag\\
\big[t^{23},[s^{14},t^{24}]\big]-\big[s^{14},[t^{23},t^{24}]\big]&=0,\notag\\
\big[s^{23},[t^{13},t^{14}]\big]-\big[[s^{23},t^{13}],t^{14}\big]&=0.
\end{align}
and
\begin{equation}\label{E:F}
\begin{split}
\big[s^{23}+s^{24},[s^{13}+s^{14},t^{34}]\big]
-\big[s^{13}&+s^{14},[s^{23}+s^{24},t^{34}]\big]\\
&=\big[[s^{23},s^{13}],t^{34}\big]+\big[[s^{24},s^{14}],t^{34}\big]\\
&=-\big[[s^{13},s^{23}],t^{34}\big]-\big[[s^{14},s^{24}],t^{34}\big].
\end{split}
\end{equation}
Collecting terms from (\ref{E:A}),(\ref{E:B}),(\ref{E:C}),(\ref{E:D}),
(\ref{E:E}),(\ref{E:F}) and replacing in (\ref{E:G}), we obtain
\begin{equation}\label{E:onze}
\begin{split}
X_i&-X_j=\\
&=(x_i^* \otimes x_j^* \otimes 1 \otimes 1)\bigg\{
-[\frac{1}{4}\big((Z_{\Omega_\l})^{123}+(Z_{\Omega_\l})^{124}\big),t^{34}]\\
&\qquad\qquad\qquad\qquad-\big[[s^{13},s^{23}],t^{34}\big]-
\big[[s^{14},s^{24}],t^{34}\big]
-[\partial_i(s^{23}+s^{24}),t^{34}]\\
&\qquad\qquad\qquad\qquad+[\partial_j(s^{13}+s^{14}),t^{34}]
+\big[s^{12},[t^{23},t^{24}]+\zom^{234}\big]\\
&\qquad\qquad\qquad\qquad+\big[s^{12},[t^{13},t^{14}]+\zom^{134}\big]
\bigg\}.
\end{split}
\end{equation}
Using the fact that $t$ is a solution to the system (\ref{E:sept}) again
we have
$$\big[s^{12},[t^{23},t^{24}]+\zom^{234}\big]+\big[s^{12},[t^{13},t^{14}]+
\zom^{134}\big]$$
\begin{equation}\label{E:quatorze}
\begin{split}
\quad&=\big[s^{12},\sum_k (x_k \otimes 1 + 1 \otimes x_k)\otimes
(-\partial_kt -[s_k \otimes 1 + 1 \otimes s_k,t])\big]\\
\quad&=-\big[s^{12},\sum_k (x_k \otimes 1 + 1 \otimes x_k)\otimes\partial_kt
\big]
-\big[s^{12},[s^{23}+s^{24},t^{34}]\big]\\
\quad&\quad-\big[s^{12},[s^{13}+s^{14},t^{34}]
\big]
\end{split}
\end{equation}
where we set $s_k=(x_k^* \otimes 1)s$. But $s$ is $\l$-equivariant, i.e for
$y \in \l$ we have
$$[s,y \otimes 1 + 1 \otimes y]=\sum_l [x_l,y]\partial_ls$$
where $[x_l,y]$ is considered as a function $D \to \C$.
Thus,
\begin{equation}\label{E:douze}
\begin{split}
-[s^{12},\sum_k (x_k \otimes 1 + 1 \otimes x_k)\otimes\partial_kt]
&=-\sum_{l,k}\partial_ls^{12}\partial_k t^{34}[x_l,x_k]\\
&=-\sum_l \partial_ls^{12}[x_l^{3}+x_l^{4},t^{34}].
\end{split}
\end{equation}
Using Jacobi identity, we can write
\begin{align}\label{E:treize}
\big[s^{12},[s^{23}+s^{24},t^{34}]\big]&=\big[[s^{12},s^{23}],t^{34}\big]
+\big[[s^{12},s^{24}],t^{34}\big],\notag\\
\big[s^{12},[s^{23}+s^{24},t^{34}]\big]&=\big[[s^{12},s^{13}],t^{34}\big]
+\big[[s^{12},s^{14}],t^{34}\big].
\end{align}
Using (\ref{E:onze}), (\ref{E:quatorze}), (\ref{E:douze}) and (\ref{E:treize})
we finally get, by (\ref{E:mcdybe})
\begin{equation*}
\begin{split}
X_i&-X_j=-(x_i^* \otimes x_j^* \otimes 1 \otimes 1)\cdot\\
&\quad\bigg\{
\big[\mathrm{Alt}\;(\overline{ds})^{123}
+[s^{13},s^{23}]+[s^{12},s^{13}]+[s^{12},s^{23}]+
\frac{1}{4}(Z_{\Omega_\l})^{123},t^{34}\big]\\
&\quad+\big[\mathrm{Alt}\;(\overline{ds})^{124}
+[s^{14},s^{24}]+[s^{12},s^{14}]+[s^{12},s^{24}]+
\frac{1}{4}(Z_{\Omega_\l})^{124},t^{34}\big]\bigg\}\\
&\quad\quad\;=0
\end{split}
\end{equation*}
\qed
\paragraph{}
Let $t=\sum t_i:\; D \to \Lambda^2\m$ be the $\l$-equivariant series
constructed by applying Lemma 4.2 succesively, starting from $t_0$.
\paragraph{}Consider the
algebraic variety
$$\mathcal{T}_\Omega=\{t \in \Lambda^2\m\;|\;CYB(t+\frac{\Omega}{2})=0\;in\;
\Lambda^3(\g/\l)\}.$$
\paragraph{}Let $x^* \in \l^*$ and consider the flow on $\Lambda^2\m$
defined by the equation 
\begin{equation}\label{E:neuf}
\frac{\partial u}{\partial \epsilon}=-(x^* \otimes 1 \otimes 1)
\bigg([u^{12},u^{13}]+[s^{12}+s^{13},u^{23}]+
\frac{1}{4}\big(CYB(\Omega)-CYB(\Omega_\l)\big)\bigg).
\end{equation}
\begin{lem}The flow (\ref{E:neuf}) preserves $\mathcal{T}_\Omega$.\end{lem}
\noindent
\textit{Proof.} Let $u \in \mathcal{T}_\Omega$. Set
$h_1=(x^* \otimes 1 \otimes 1)\big([s^{12}+s^{13},u^{23}]\big)$,
$$h_2=(x^* \otimes 1 \otimes 1)
\bigg([u^{12},u^{13}]+\frac{1}{4}\big(Z_{\Omega}-Z_{\Omega_\l}\big)\bigg).$$
Note that  $h_1 \in \Lambda^2\m$ by condition i) and that
$h_2 \in \Lambda^2\m$ since $u \in \Lambda^2\m$ and since by ii),
$$Z_{\Omega}-Z_{\Omega_\l}\in (\m \otimes \g\otimes \g) \oplus (\l \otimes
\m \otimes \m).$$
It thus remains to check that
the vector field defined by (\ref{E:neuf}) is tangent to
$\mathcal{T}_\Omega$, i.e that $CYB(u,h_1+h_2)\in \mathrm{Alt}\;(\l
\otimes \g \otimes \g)$,
where we use the notation
$$CYB(a,b)=[a^{12},b^{13}]+[a^{13},b^{23}]+[a^{12},b^{23}]+
[b^{12},a^{13}]+[b^{13},a^{23}]+[b^{12},a^{23}].$$
But $$CYB(u,h_1)=\mathrm{ad}\;\big((x^*\otimes 1)s\big) CYB(u) \in
\mathrm{Alt}\;(\l\otimes \g \otimes \g),$$
and $CYB(u,h_2)\in \mathrm{Alt}\;(\l\otimes \g \otimes \g)$ by \cite{S},
Lemma 3 (note that the commutativity of $\l$, assumed in \cite{S},
is not used in the proof of Lemma 3).\qed

\begin{cor}\label{L:1} The map $t : D \to \Lambda^2\m$ takes values in 
$\mathcal{T}_\Omega$.\end{cor}
\noindent
\textit{Proof.} Note that $t(0) \in \mathcal{T}_\Omega$ by assumption, and
that for any $x^* \in D$ the function $u(\epsilon)=t(\epsilon x^*)$ on the
formal disc satisfies 
(\ref{E:neuf}) by construction. Hence $t$ takes values in $\mathcal{T}_\Omega$.
\qed

\paragraph{}We now conclude the proof of Theorem 3 by showing that 
$r=r^\l_{AM}+\frac{\Omega_\m}{2}+t$ is a dynamical r-matrix. Setting
$s^\l_{AM}=r^\l_{AM}-\frac{\Omega_\l}{2}$ we have
\begin{equation*}
\begin{split}
\mathrm{Alt}&\;(\overline{dr})+CYB(r)\\
&=\mathrm{Alt}\;(\overline{ds^\l_{AM}})
+\mathrm{Alt}\;(\overline{dt}) + CYB(s^{\l}_{AM}+t) + \frac{1}{4}CYB(\Omega)\\
&=\mathrm{Alt}\;(\overline{ds^\l_{AM}})
+\mathrm{Alt}\;(\overline{dt}) + CYB(s^{\l}_{AM})+CYB(t)+CYB(s^{\l}_{AM},t)\\
&\quad+\frac{1}{4}CYB(\Omega).
\end{split}
\end{equation*}
Using the CDYBE for $r^\l_{AM}$ we see that $r$ is a dynamical r-matrix
if and only if
\begin{equation}\label{E:huit}
\mathrm{Alt}\;(\overline{dt})+CYB(t)+ CYB(s^{\l}_{AM},t)+
\frac{1}{4}\big(CYB(\Omega)-CYB(\Omega_\l)\big)=0.
\end{equation}
Since $t$ takes values in $\mathcal{T}_\Omega$, (\ref{E:huit}) is equivalent
to the system
$$\frac{\partial t}{\partial x_i^*}=-(x_i^* \otimes 1 \otimes 1)
\bigg(CYB(t)+ CYB(s^{\l}_{AM},t)+
\frac{1}{4}\big(CYB(\Omega)-CYB(\Omega_\l)\big)\bigg)$$
for $i=1,\ldots r$. It is easy to see from conditions i) and ii) that this
last system is itself equivalent to the collection of systems
(\ref{E:sept}) for all $l \in \N$.\qed

\section{Appendix. Generalized Alekseev-Meinrenken\\ dynamical r-matrices}
\paragraph{}In this appendix we give a generalization of the dynamical
r-matrix
$r^\l_{AM}$.
\paragraph{}Let $\g$ be a finite-dimensional complex Lie algebra and $B:
\g \to \g$ an automorphism of order $n$. Then $\g=\bigoplus_{j \in \Z/n\Z}
\g_j$ where $\g_j=\mathrm{Ker}\;(B-e^{\frac{2i\pi j}{n}})$. Set $\l=\g_0$.
Then $\g_0$ acts on $\g_j$ for all $j$.\\
\hbox to1em{\hfill}Assume that $\g$ carries a nondegenerate invariant
form $(\,,\,)$, which is stable under $B$. Set $\Omega=(\,,\,)^{-1}
\in (S^2\g)^\g$. We will identify $\g$ with $\g^*$ and
$\l$ with $\l^*$ using $(\,,\,)$.
\paragraph{}Let $D$ be the formal neighborhood of zero in $\l^*\simeq \l$.
Consider the function $\hat{\rho}:D \to \mathrm{End}(\g)$ such that
$\hat{\rho}(A)_{|\g_i}=f_i(\mathrm{ad}\;A),$
with
\begin{align*}
f_0(s)&=\frac{1}{s}-\frac{1}{2}\mathrm{cotanh}\;(\frac{1}{2}s),\\
f_j(s)&=-\frac{1}{2}\mathrm{cotanh}\;(\frac{1}{2}(s+\frac{2i\pi j}{n})),
\qquad j \neq 0.
\end{align*}
The element $\hat{\rho}$ defines a map
$\rho: D \to \Lambda^2\g$. Let us set $r_B=\frac{\Omega}{2}+\rho$.
\begin{theoA} The map $r_B$ is a dynamical r-matrix.\end{theoA}
\paragraph{Remark 7.} If $B=1$ then $r_B$ is equal to the Alekseev-Meinrenken
dynamical r-matrix $r^\l_{AM}$.
\paragraph{}The rest of this appendix is devoted to the proof of Theorem A.1.
We start by recalling the following result from \cite{EV}. Let $\l$ be a
reductive Lie algebra with Cartan subalgebra $\h$, $\g$ any finite-dimensional
Lie algebra containing $\l$ and let $\Omega \in (S^2\g)^\g$. The projection
$\l \to \h$ defines an embedding $\h^* \to \l^*$. 
Let $\Delta$ be the root system of $\l$ and $\l_\alpha$
the weight subspace corresponding to $\alpha \in \Delta$. 
Choose an nondegenerate invariant inner product on $\l$. 
Let us fix $e_\alpha \in \l_\alpha$ for all $\alpha \in \Delta$ such
that $(e_{-\alpha},e_\alpha)=1$. Define a function $\rho_0: D \to \Lambda^2\l
 \subset \Lambda^2\g$ by
$$\rho_0(\lambda)=\sum_{\alpha>0} 
\frac{e_{\alpha}\otimes e_{-\alpha}-e_{\alpha}\otimes
e_{-\alpha}}{(\alpha,\lambda)}.$$
It is clear that $\rho_0$ does not depend on the choice 
of the inner product. 

Let
$r: \l^* \to \g \otimes \g$ be an $\l$-equivariant meromorphic function
satisfying the quasi-unitarity condition $r+r^{21}=\Omega$.
\begin{theoA}[\cite{EV}, Theorem 3.14] The map $r$ is a classical dynamical
r-matrix if and only if $r_{|\h^*} + \rho_0$ is a classical dynamical r-matrix
for $\h$.\end{theoA}
\noindent
\textit{Proof.} This is proved in \cite{EV} under the assumption that $\g$
is simple and $\h \subset \g$ is a Cartan subalgebra. However, this assumption
is not used in the proof and the result is valid in general.\qed
\begin{propA} Theorem A1 is valid if $\l$ is reductive, $\g=\l_1 \oplus
\cdots \oplus \l_n$ with $\l_i=\l$, and $B$ is the cyclic permutation
automorphism $B: \l_i \stackrel{\sim}{\to} \l_{i+1\;\mathrm{mod}\;n}$.
\end{propA}
\noindent
\textit{Proof.} Let $(x_i)_{i\in I}$
be an orthonormal basis of $\h \subset \l$. For $i=1,\ldots,n$ we will write
$e_{\alpha}^{(i)}$ for the element of $\l_i\subset \g$ corresponding to
$e_\alpha$. With this notation, we have
$B(e_{\alpha}^{(i)})=e_{\alpha}^{(i+1)}$ and
$$\g_j=\{g_1^{(1)}\oplus\cdots\oplus g_n^{(n)};|\;g_{k+1}=
e^{\frac{-2i\pi j}{n}}g_{k}\}.$$
Finally, let $\mathfrak{k}\subset \bigoplus_i \h^{(i)}$ be the orthogonal
complement to $\h$. Note that $1-B$ restricts to an invertible operator on
$\mathfrak{k}$. A direct computation shows that
\begin{equation*}
\begin{split}
r_{B|\h^*}&+\rho_0\\
&=\frac{\Omega}{2}+\sum_{\alpha >0,i}
e_{\alpha}^{(i)} \wedge \big(-\frac{1}{2}\big(
\frac{1+Be^{(\alpha,\lambda)}}{1-Be^{(\alpha,\lambda)}}\big)e_{-\alpha}^{(i)}
\big)+\sum_i y_i \otimes \big(-\frac{1+B}{2(1-B)}y_i\big)\\
&=\sum_i x_i \otimes x_i+\sum_{\alpha>0,i} e^{(i)}_{-\alpha} \otimes
e_{\alpha}^{(i)} -\sum_{\alpha >0,i}\sum_{l \geq 1}e^{l(\alpha,\lambda)}
e_{\alpha}^{(i)} \wedge e^{(i+l)}_{-\alpha}\\
&\quad +\frac{1}{2}\sum_i \frac{B+1}{B-1}
y_i \otimes y_i
\end{split}
\end{equation*}
where $(y_i)_{i \in J}$ is an orthonormal basis of $\mathfrak{k}$.
By \cite{S}, Theorem 4 this expression is a dymamical r-matrix. Hence, by
Theorem A.2, $r_B$ is a dynamical r-matrix.\qed

\paragraph{}Define a map $W:D \to \Lambda^3\g$ by
$$W(A)=\mathrm{Alt}\;(\overline{d\rho(A)})+CYB(\rho(A))+\frac{1}{4}Z$$
where $Z=CYB(\Omega)$. For any $i,j \in \Z/n\Z$, $X \in \g_i$, 
$Y\in \g_j$ and $A \in \l$, consider the expression
$$K_{ij}(A,X,Y)=\big(1 \otimes X \otimes Y,W(A)\big) \in \g_{i+j}.$$

\begin{lemA} The expression $K_{ij}(A,X,Y)$ is given by a universal Lie
series in $A,X$ and $Y$.\end{lemA}
\noindent
\textit{Proof.} Straightforward.\qed
\paragraph{}Moreover, from Proposition A.1 we deduce the following result.
\begin{propA} We have $K_{ij}(A,X,Y)=0$ for all $A,X,Y$ if
$\l=\mathfrak{gl}_n(\C)$, $\g=\l_1 \oplus
\cdots \oplus \l_n$ with $\l_k=\l$, and $B$ is the cyclic permutation
automorphism.\end{propA}
\paragraph{}Finally, we recall the following standard fact.
\begin{lemA} Let $P(X_1,\ldots,X_n)$ be a Lie polynomial which
vanishes identically for all $X_1,\ldots X_n \in \mathfrak{gl}_k(\C)$ for all
$k \in \N$. Then $P(X_1,\ldots,X_n)=0$.\end{lemA}
\noindent
\textit{Proof.} Let $F_n$ be the free Lie algebra in $n$ generators and let
$U_n$ be its enveloping algebra (the free associative algebra). Let $d$ be the
degree of $P$ and let $I$ be the ideal in $U_n$
generated by elements of degree at least
$d+1$. Then $U_n/I$ is a finite-dimensional algebra. Let
$\sigma : U_n/I\to \mathfrak{gl}(U_n/I)$ be the left regular representation.
Then $\sigma\big(P(X_1,\ldots,X_n)\big)=0$. Hence $P(X_1,\ldots,X_n)=0$.\qed
\paragraph{}Now, let us write $K_{ij}=\sum_k K_{ij}^{(k)}$ where $K_{ij}^{(k)}$
is the homogeneous component of degree $k$. By Proposition A.2, 
$K_{ij}^{(k)}(A,X,Y)=0$ whenever $A,X,Y \in \mathfrak{gl}_m(\C)$ for some $m
\in \C$. Hence $K_{ij}^{(k)}=0$ by Lemma A.2. Thus $W(A)=0$ for all $A \in \l$.
Theorem A.1 is proved.
\paragraph{Examples.} Let $\g$ be a simple complex Lie algebra
and $\l \subset \g$
a semisimple subalgebra with same rank as $\g$. Such pairs are
classified in \cite{BDS}. Let $Q_\l$ and $Q_\g$ be the root lattices of $\l$
and $\g$ respectively and set $\Gamma=Q_{\g}/Q_{\l}$. It follows from
\cite{BDS} that $\Gamma$ is one of the groups $\Z/2\Z$, $\Z/3\Z$ or $\Z/5\Z$
(the case $\Gamma=\Z/2\Z$ corresponds to symmetric spaces). Let
$\chi$ be a nontrivial character $\chi$ of $\Gamma$. Then $\chi$
 gives rise to an automorphism $B_\chi$ of $\g$ whose set of fixed points
is $\l$, defined by
$$B_{\chi|\g_\alpha}=\chi(\alpha)Id,$$
where $\g_\alpha$ is the root space of weight $\alpha$. Let $\Omega \in
(S^2\g)^\g$ be a Casimir element ($\Omega \neq 0$).
Let $r: D \to \g \otimes \g$ be a dynamical r-matrix such that
$r+r^{21}=\Omega$. It follows from Theorem A.2 and \cite{EV} Theorem 3.1, that,
up to gauge-equivalence,
$$r_{|\h^*}+\rho_0=\frac{\Omega}{2}+\sum_{\alpha \in \Delta_\g}
\frac{1}{2}\mathrm{cotanh}(\frac{1}{2}(\alpha,\lambda-\nu))e_{\alpha}
\otimes e_{-\alpha},$$
for some $\nu \in \h^*$, where $\Delta_\g$ is the root system of $\g$. But then
$r$ is regular at $\lambda=0$ if and only if $(\alpha,\nu)=0$ 
modulo $2\pi i\Z$ for all
$\alpha \in \Delta_\l$ and $(\alpha,\nu) \neq 0$ 
modulo $2\pi i\Z$ for all $\alpha\in
\Delta_{\g}\backslash\Delta_\l$, where $\Delta_\l\subset \Delta_\g$ is the
root system of $\l$. Such $\nu$ defines a nontrivial character $\chi$ of
$\Gamma$, and it follows from \cite{EV}, Section 3.8 that $r$ is
gauge-equivalent to the generalized
Alekseev-Meinrenken dynamical r-matrix $r_{B_\chi}$.
Hence the moduli space $\mathcal{M}(\g,\l,\Omega)$
consists of $|\Gamma|-1$ points.

\paragraph{Acknowledgments.} The authors are grateful to P. Xu for 
useful discussions. The first author was supported by the NSF grant 
DMS-9700477. P.E performed this research as a CMI prize fellow. O.S conducted
this research partially for the Clay Mathematics Institute, and thanks the MIT
mathematics department for hospitality.

\end{document}